\documentclass{amsart}

\usepackage{amsmath}
\usepackage{amsfonts}
\usepackage{amssymb}

\numberwithin{equation}{section}

\newtheorem{theorem}[subsection]{Theorem}
\newtheorem{proposition}[subsection]{Proposition}

\newtheorem{lemma}[subsection]{Lemma}

\newtheorem{conjecture}[subsection]{Conjecture}

\newtheorem{question}[subsection]{Question}

\newcommand{\wh}{\widehat}

\title[The A(G)-norm on compact vector spaces]
 {The $\ell^1$-norm of the Fourier transform on compact vector spaces}

\author{T. Sanders}
\address{Department of Pure Mathematics and Mathematical Statistics\\
University of Cambridge\\
Wilberforce Road\\
Cambridge CB3 0WA\\
England } \email{t.sanders@dpmms.cam.ac.uk}

\begin{document}

\begin{abstract}
Suppose that $G$ is a compact Abelian group. If $A \subset G$ then
how small can $\|\chi_A\|_{A(G)}$ be? In general there is no
non-trivial lower bound.

In \cite{BJGSVK} Green and Konyagin show that if $\widehat{G}$ has
sparse small subgroup structure and $A$ has density $\alpha$ with
$\alpha(1-\alpha) \gg 1$ then $\|\chi_A\|_{A(G)}$ does admit a
non-trivial lower bound.

In this paper we address the complementary case of groups with
duals having rich small subgroup structure, specifically the case
when $G$ is a compact vector space over $\mathbb{F}_2$. The
results themselves are rather technical to state but the following
consequence captures their essence: If $A \subset \mathbb{F}_2^n$
is a set of density as close to $\frac{1}{3}$ as possible then we
show that $\|\chi_A\|_{A(\mathbb{F}_2^n)} \gg \log n$.

We include a number of examples and conjectures which suggest that
what we have shown is very far from a complete picture.
\end{abstract}

\maketitle

\section{Notation and introduction}\label{intro}

We use the Fourier transform on compact Abelian groups, the basics
of which may be found in Chapter 1 of Rudin \cite{WR}; we take a
moment to standardize our notation.

Suppose that $G$ is a compact Abelian group. Write $\widehat{G}$
for the dual group, that is the discrete Abelian group of
continuous homomorphisms $\gamma:G \rightarrow S^1$, where
$S^1:=\{z \in \mathbb{C}:|z|=1\}$. Although the natural group
operation on $\wh{G}$ corresponds to pointwise multiplication of
characters we shall denote it by '$+$' in alignment with
contemporary work. $G$ may be endowed with Haar measure $\mu_G$
normalised so that $\mu_G(G)=1$ and as a consequence we may define
the Fourier transform $\widehat{.}:L^1(G) \rightarrow
\ell^\infty(\widehat{G})$ which takes $f \in L^1(G)$ to
\begin{equation*}
\widehat{f}: \widehat{G} \rightarrow \mathbb{C}; \gamma \mapsto
\int_{x \in G}{f(x)\overline{\gamma(x)}d\mu_G(x)}.
\end{equation*}
\noindent We write
\begin{equation*}
A(G):=\{f \in L^1(G): \|\widehat{f}\|_1 < \infty\},
\end{equation*}
and define a norm on $A(G)$ by $\|f\|_{A(G)}:=\|\widehat{f}\|_1$.

The following proposition is easy to prove and has been known at
least since the 1960s.
\begin{proposition}\label{onethirdold}
Suppose that $G$ is a compact vector space over $\mathbb{F}_2$ and
$A \subset G$ has density (i.e. measure with respect to the
probability measure $\mu_G$) $1/3$. Then $\chi_A \not \in A(G)$.
\end{proposition}
This proposition only has content for infinite $G$; in this paper
we prove the following quantitative version which has content for
finite $G$ too.
\begin{theorem}\label{onethird}
Suppose that $G$ is a compact vector space over $\mathbb{F}_2$ and
$A \subset G$ has density $\alpha$ with $|\alpha - 1/3| \leq
\epsilon$. Then
\begin{equation*}
\|\chi_A\|_{A(G)} \gg \log \log \epsilon^{-1}.
\end{equation*}
\end{theorem}
In view of Proposition \ref{onethirdold} one could restrict
attention in Theorem \ref{onethird} to the case when $G$ is finite
rendering the compactness requirement irrelevant.

There is a simple construction, the details of which are given in
\S\ref{egs}, of a set $A$ satisfying the hypotheses of the theorem
with $\|\chi_A\|_{A(G)} \ll \log \epsilon ^{-1}$; we conjecture
that this construction represents the true state of affairs
\emph{viz}.
\begin{conjecture}
Suppose that $G$ is a compact vector space over $\mathbb{F}_2$ and
$A \subset G$ has density $\alpha$ with $|\alpha - 1/3| \leq
\epsilon$. Then
\begin{equation*}
\|\chi_A\|_{A(G)} \gg \log \epsilon^{-1}.
\end{equation*}
\end{conjecture}

The paper splits into five further sections. Although the result
of the introduction captures the spirit of this work, in
\S\ref{context} we explain our results in a more general and
natural context. \S\ref{egs} then provides some examples which
complement our results and are worth bearing in mind when
following the proof. \S\ref{finfourier} is the central iterative
argument; in this section we essentially prove a result with the
conclusion of Theorem \ref{onethird} but with a more cumbersome
hypothesis on $A$. \S\ref{finphysical} then provides some physical
space estimates to show that sets of density close to 1/3 satisfy
this hypothesis. The final section, \S\ref{finrem}, combines the
work of the previous two sections to prove a result which implies
the main result of the paper and discusses the limitations of our
methods.

\section{The problem in context}\label{context}

Proposition \ref{onethirdold} is, in fact, a special case of the
following equally simple but more natural result essentially due
to Cohen \cite{PJC}; a proof and more detailed discussion is
contained in \cite{TS.ALVCT}.

\begin{proposition}\label{Cohen}
Suppose that $G$ is a compact Abelian group. Suppose that $A
\subset G$ has density $\alpha$ and for all finite $V \leq \wh{G}$
we have $\{\alpha|V|\}(1-\{\alpha|V|\})>0$ where $\{\alpha|V|\}$
denotes the fractional part of $\alpha|V|$. Then $\chi_A \not \in
A(G)$.
\end{proposition}
We are interested in quantitative versions of this proposition.
Specifically we ask the following question.
\begin{question}\label{main question}
Suppose that $G$ is a compact Abelian group. Suppose that $A
\subset G$ has density $\alpha$ and for all finite $V \leq \wh{G}$
with $|V| \leq M$ we have $\{\alpha|V|\}(1-\{\alpha|V|\})\gg 1$.
Then how small can $\|\chi_A\|_{A(G)}$ be in terms of $M$?
\end{question}
Green and Konyagin made progress on this question in \cite{BJGSVK}
where they addressed the case when $\wh{G}$ has sparse small
subgroup structure. They proved the following result.
\begin{theorem}\label{bjgsk}
Suppose that $G$ is a compact Abelian group and the only subgroup
$V \leq \wh{G}$ with $|V|\leq M$ is the trivial group. Suppose
that $A \subset G$ has density $\alpha$ and for all finite $V \leq
\wh{G}$ with $|V| \leq M$ we have
$\{\alpha|V|\}(1-\{\alpha|V|\})\gg 1$. Then
\begin{equation*}
\|\chi_A\|_{A(G)} \gg \left(\frac{\log M}{\log \log
M}\right)^{\frac{1}{3}}.
\end{equation*}
\end{theorem}
The density condition here collapses to $\alpha(1-\alpha)\gg 1$
but we retain the more complex form for illustrative purposes. In
this paper we address the complementary case when $\wh{G}$ has
very rich small subgroup structure, specifically when $G$ is a
compact vector space over $\mathbb{F}_2$. We prove the following
result.
\begin{theorem}\label{maintheorem}
Suppose that $G$ is a compact vector space over $\mathbb{F}_2$.
Suppose that $A \subset G$ has density $\alpha$ and for all finite
$V \leq \wh{G}$ with $|V| \leq M$ we have
$\{\alpha|V|\}(1-\{\alpha|V|\})\gg 1$. Then
\begin{equation*}
\|\chi_A\|_{A(G)} \gg \log \log M.
\end{equation*}
\end{theorem}
The ideas of this paper and \cite{BJGSVK} can be combined to
address the original question more fully. In \cite{TS.ALVCT} we show
the following result.
\begin{theorem}
Suppose that $G$ is a compact Abelian group. Suppose that $A
\subset G$ has density $\alpha$ and for all finite $V \leq \wh{G}$
with $|V| \leq M$ we have $\{\alpha|V|\}(1-\{\alpha|V|\})\gg 1$.
Then
\begin{equation*}
\|\chi_A\|_{A(G)} \gg \log \log \log M.
\end{equation*}
\end{theorem}
The examples of \S\ref{egs} show that the bound in Theorem
\ref{maintheorem} cannot be better than a constant multiple of
$\log M$, and there are easy examples mentioned by Green and
Konyagin in \cite{BJGSVK} to show that the bound in Theorem
\ref{bjgsk} cannot be better than a constant multiple of $\log M$.
In the absence of any better examples we make the following
conjecture.
\begin{conjecture}
Suppose that $G$ is a compact Abelian group. Suppose that $A
\subset G$ has density $\alpha$ and for all finite $V \leq \wh{G}$
with $|V| \leq M$ we have $\{\alpha|V|\}(1-\{\alpha|V|\})\gg 1$.
Then
\begin{equation*}
\|\chi_A\|_{A(G)} \gg \log M.
\end{equation*}
\end{conjecture}

The various statements in the introduction follow easily from
those of this section once one recalls that if $G$ is a compact
vector space over $\mathbb{F}_2$ and $V \leq \wh{G}$ is finite
then $|V|$ is a power of 2.

\section{Sets with small $A(G)$-norm}\label{egs}

Throughout this section $G$ is a compact vector space over
$\mathbb{F}_2$.

We address the question of how to construct subsets of $G$ of a
prescribed density whose characteristic function has small
$A(G)$-norm. This complements the main results of the paper.

Cosets are the simplest example of sets whose characteristic
function has small $A(G)$-norm: Recall that if $V \leq
\widehat{G}$ then we write $V^\perp$ for the \emph{annihilator} of
$V$ i.e.
\begin{equation*}
V^\perp:=\{x \in G: \gamma(x)=1 \textrm{ for all } \gamma \in V\}.
\end{equation*}
If $V \leq \widehat{G}$ is finite and $A=x+V^\perp$ then a simple
calculation gives
\begin{equation*}
\widehat{\chi_A}(\gamma)=\begin{cases} \gamma(x)|V|^{-1} &
\textrm{ if }\gamma \in V\\ 0 & \textrm{ otherwise.}
\end{cases}
\end{equation*}
It follows that $\|\chi_A\|_{A(G)}=1$. A coset has density
$2^{-d}$ for some integer $d$; to produce a set with a density not
of this form we take unions of cosets.

Suppose that we are given $\alpha \in [0,1]$ a terminating binary
number. Write
\begin{equation*}
\alpha=\sum_{i=1}^k{2^{-d_i}},
\end{equation*}
where the $d_i$ are strictly increasing. If we can find a sequence
of disjoint cosets $A_1,...,A_k$ such that $\mu_G(A_i)=2^{-d_i}$,
then their union $A:=\bigcup_{i=1}^k{A_i}$ has
\begin{equation}\label{unionbound}
\|\chi_A\|_{A(G)} = \|\sum_{i=1}^k{\chi_{A_i}}\|_{A(G)} \leq
\sum_{i=1}^k\|\chi_{A_i}\|_{A(G)} = k
\end{equation}
by the triangle inequality, and density
\begin{equation*}
\mu_G(A)=\sum_{i=1}^k{\mu_G(A_i)}= \sum_{i=1}^k{2^{-d_i}} = \alpha
\end{equation*}
since the elements of the union are disjoint. To produce such
cosets we take $\{0_{\wh{G}}\} = \Lambda_0 < \Lambda_1 < ... <
\Lambda_k \leq \wh{G}$ a nested sequence of subspaces with $\dim
\Lambda_i=d_i$. Choose a sequence of vectors $\{\gamma_i:1 \leq i
\leq k\}$ such that $\gamma_i \in \Lambda_i \setminus
\Lambda_{i-1}$ for $1 \leq i \leq k$. It is easy to see that this
sequence must be linearly independent so we may take a sequence
$\{x_i:1 \leq i \leq k-1\}$ such that
\begin{equation}\label{xis}
\gamma_j(x_i)=\begin{cases}  1 & \textrm{ if } j \neq i\\ -1 &
\textrm{ if } j = i
\end{cases}
\end{equation}
for all $1 \leq i \leq k-1$. Put
\begin{equation*}
A_i=x_1+...+x_{i-1}+\Lambda_i^\perp.
\end{equation*}
First we note that $\mu_G(A_i)=2^{-d_i}$ and second that the sets
$A_i$ are pairwise disjoint: Suppose $j>i$ and $x \in A_j$ then
$x=x_1+...+x_{j-1}+x'$ where $x' \in \Lambda_j^\perp$ so that
\begin{equation*}
\gamma_i(x)=\gamma_i(x_1)...\gamma_i(x_{j-1}).\gamma_i(x').
\end{equation*}
$i<j$ so $\Lambda_i < \Lambda_j$ from which it follows that
$\gamma_i(x')=1$. Consequently $\gamma_i(x) = \gamma_i(x_i) = -1$
by (\ref{xis}). However if $x \in A_i$ then by a similar
calculation $\gamma_i(x)=1$.

It follows that $A_1,...,A_k$ are disjoint cosets of the
appropriate size and hence their union, $A:=\bigcup_{i=1}^k{A_i}$,
has density $\alpha$ and $\|\chi_A\|_{A(G)} \leq k$.

The first density we apply this general construction to is
\begin{equation*}
\alpha= \frac{1}{4}+\frac{1}{16} +... + \frac{1}{4^k}.
\end{equation*}
The set $A$ we produce has density $\alpha$ and the following two
properties.
\begin{enumerate}
\item  $A$ satisfies the hypotheses of Theorem \ref{maintheorem}
with $M=4^k-1$: If $V \leq \wh{G}$ and $|V| \leq M$ then $|V|=2^d$
for some $d < k$ and
\begin{equation*}
\frac{2}{3} \geq \sum_{i=\lfloor d/2\rfloor +1}^k{2^d.4^{-i}} =
\{\alpha 2^d\} \geq \frac{2^d}{4^{\lfloor d/2\rfloor +1}} \geq
\frac{1}{4},
\end{equation*}
and hence $\{\alpha|V|\}(1-\{\alpha|V|\}) \geq 1/12$.
\item $\|\chi_A\|_{A(G)} \asymp k$: $\|\chi_A\|_{A(G)} \leq k$
follows by construction; $\|\chi_A\|_{A(G)} \gg k$ is slightly
more involved:
\begin{equation*}
\wh{\chi_{A_i}}(\gamma)=\begin{cases}
4^{-i}\gamma(x_1)...\gamma(x_{i-1}) & \textrm{ if
} \gamma \in \Lambda_i\\
0 & \textrm{ otherwise.}
\end{cases}
\end{equation*}
Hence we can bound $|\wh{\chi_A}|$ from below using the linearity
of the Fourier transform.
\begin{equation*}
|\wh{\chi_A}(\gamma)| \geq 4^{-i}- \sum_{j=i+1}^k{4^{-j}} \geq
\frac{2}{3}.4^{-i} \textrm{ if } \gamma \in \Lambda_i\setminus
\Lambda_{i-1},
\end{equation*}
and so
\begin{equation*}
\|\chi_A\|_{A(G)} \geq
\sum_{i=1}^k{\frac{2}{3}4^{-i}.|\Lambda_i\setminus \Lambda_{i-1}|}
\geq \sum_{i=1}^k{\frac{2}{3}4^{-i}.\frac{3}{4}4^i} = \frac{k}{2}.
\end{equation*}
\end{enumerate}
The conclusion of Theorem \ref{maintheorem} is that
$\|\chi_A\|_{A(G)} \gg \log k$ which should be compared with the
fact that actually $\|\chi_A\|_{A(G)} \asymp k$.

\section{An iteration argument in Fourier space}\label{finfourier}

Throughout this section $G$ is a compact vector space over
$\mathbb{F}_2$ and $A \subset G$ has density $\alpha$.

\subsection{A trivial lower bound}

Suppose that $\alpha>0$. It is natural to try to bound
$\|\chi_A\|_{A(G)}$ by a combination of H\"{o}lder's inequality
and Plancherel's theorem:
\begin{equation}\label{ht}
\|\chi_A\|_{A(G)}\|\widehat{\chi_A}\|_\infty \geq
\|\widehat{\chi_A}\|_2^2=\|\chi_A\|_2^2;
\end{equation}
non-negativity of $\chi_A$ means that
$\widehat{\chi_A}(0_{\widehat{G}})=\|\chi_A\|_1$ so
\begin{equation}\label{calcn}
\|\chi_A\|_1 \geq \|\widehat{\chi_A}\|_\infty \geq
\widehat{\chi_A}(0_{\widehat{G}})=\|\chi_A\|_1 \Rightarrow
\|\widehat{\chi_A}\|_\infty=\|\chi_A\|_1.
\end{equation}
$\chi_A \equiv \chi_A^2$ so $\|\chi_A\|_2^2=\|\chi_A\|_1=\alpha$,
which is positive, and hence (\ref{ht}) tells us that
\begin{equation}
\label{trivialbound} \|\chi_A\|_{A(G)} \geq 1.
\end{equation}
Taking $A=G$ shows that we can do no better than this, so for
$\|\chi_A\|_{A(G)}$ to be large we require some sort of (further)
restriction on what $A$ may be. The restriction we take evolves
out of modifying the idea above so that it becomes a step in an
iteration.

\subsection{A weak iteration lemma}

A weakness in the above deduction is that we have no good upper
bound for $\|\widehat{\chi_A}\|_\infty$. In fact, as we saw,
$\|\widehat{\chi_A}\|_\infty$ is necessarily large because
$\widehat{\chi_A}$ is large at the trivial character, however we
know nothing about how large $\widehat{\chi_A}$ is at any other
character, a fact which we shall now exploit.

Write $f$ for the \emph{balanced function} of $\chi_A$ i.e.
$f=\chi_A - \alpha$. Then
\begin{equation*}
\widehat{f}(\gamma)=\begin{cases} 0 &\textrm{ if } \gamma =
0_{\widehat{G}}\\ \widehat{\chi_A}(\gamma) & \textrm{ otherwise.}
\end{cases}
\end{equation*}
Applying H\"{o}lder's inequality and Plancherel's theorem in the
same way as before we have
\begin{equation}\label{h2}
\|\chi_A\|_{A(G)}\|\widehat{f}\|_\infty \geq \langle
\widehat{\chi_A},\widehat{f}\rangle=\langle \chi_A, f
\rangle=\alpha-\alpha^2.
\end{equation}
Now, fix $\epsilon>0$ to be optimized later. If $\alpha$ is
bounded away from 0 and 1 by an absolute constant then either
$\|\chi_A\|_{A(G)} \gg \epsilon^{-1}$ or $\|\widehat{f}\|_\infty
\gg \epsilon$. In the former case we are done (since
$\|\chi_A\|_{A(G)}$ is large) and in the latter we have a
non-trivial character at which $\widehat{\chi_A}$ is large; we
should like to start building up a collection of such characters.

Suppose that $\Gamma \subset \widehat{G}$ is a collection of
characters on which we know $\widehat{\chi_A}$ has large
$\ell^1$-mass. We want to produce a superset $\Gamma'$ of $\Gamma$
by adding some more characters which support a significant
$\ell^1$-mass of $\widehat{\chi_A}$. To find characters outside
$\Gamma$ on which $\widehat{\chi_A}$ has large $\ell^1$-mass we
might replace $f$ with a function $f_{\Gamma}$ (in analogy with
the earlier replacement of $\chi_A$ by $f$) defined by inversion:
\begin{equation}\label{tacomp}
\widehat{f_{\Gamma}}(\gamma) = \begin{cases} 0 & \textrm{ if }
\gamma \in \Gamma\\ \widehat{\chi_A}(\gamma) & \textrm{
otherwise.}
\end{cases}
\end{equation}
The problem with this is that for general $\Gamma$ we can say very
little about $f_{\Gamma}$. If $V \leq \widehat{G}$, however, then
$f_V$ has a particularly simple form:
\begin{equation*}
f_{V}=\sum_{\gamma \in
\widehat{G}}{\widehat{\chi_A}(\gamma)(1-\widehat{\mu_{V^\perp}}(\gamma))\gamma}=\chi_A
\ast ( \delta - \mu_{V^\perp}) = \chi_A - \chi_A \ast
\mu_{V^\perp} \textrm{ a.e.}
\end{equation*}
Now suppose that $\Gamma=V \leq \widehat{G}$. We want to try to
add characters to $V$ to get a superspace $V' \leq \widehat{G}$
with
\begin{equation*}
\sum_{\gamma \in V'}{|\widehat{\chi_A}(\gamma)|} \textrm{
significantly larger than } \sum_{\gamma \in
V}{|\widehat{\chi_A}(\gamma)|}.
\end{equation*}
We can use the idea in (\ref{h2}) to do this; replace $f$ by
$f_{V}$ in that argument:
\begin{equation}\label{h3}
\|\chi_A\|_{A(G)}\|\widehat{f_{V}}\|_\infty \geq \langle
\widehat{\chi_A},\widehat{f_{V}}\rangle=\langle \chi_A, f_{V}
\rangle.
\end{equation}
Before, an easy calculation gave us $\langle
\chi_A,f\rangle=\alpha(1-\alpha)$. To compute $\langle \chi_A,f_V
\rangle$ we have a slightly more involved calculation.
\begin{lemma}\label{t_Alem}
\begin{equation}
\label{t_A} \|f_{V}\|_1 = 2 \langle \chi_A,f_{V} \rangle.
\end{equation}
\end{lemma}
\begin{proof}
$\mu_{V^\perp}$ is a probability measure so $0 \leq \chi_A \ast
\mu_{V^\perp}(x) \leq 1$, hence $f_{V}(x)\leq 0$ for almost all $x
\not \in A$ and $f_{V}(x)\geq 0$ for almost all $x \in A$;
consequently
\begin{equation*}
\|f_{V}\|_1 = \int{\chi_Af_{V}d\mu_G} +
\int{(1-\chi_A)(-f_{V})d\mu_G} = 2\langle \chi_A,f_{V} \rangle -
\int{f_{V}d\mu_G}.
\end{equation*}
But $\int{f_{V}d\mu_G}=0$ since $\int{\chi_Ad\mu_G}=\int{\chi_A
\ast \mu_{V^\perp}d\mu_G}$, so we are done.
\end{proof}
It follows that
\begin{equation}
\|\chi_A\|_{A(G)}\|\widehat{f_{V}}\|_\infty
\geq\frac{\|f_{V}\|_1}{2}.
\end{equation}
So either $\|\chi_A\|_{A(G)} \geq \epsilon^{-1}$ or there is a
character $\gamma$ such that $|\widehat{f_V}(\gamma)| \geq
\epsilon \|f_V\|_1/2$. By construction of $f_V$ we have
$\widehat{f_V}(\gamma')=0$ if $\gamma' \in V$ so that $\gamma \not
\in V$ -- $\gamma$ is a genuinely new character. We let $V'$ be
the space generated by $\gamma$ and $V$ and have our first
iteration lemma:
\begin{lemma}
\emph{(Weak iteration lemma)} Suppose that $G$ is a compact vector
space over $\mathbb{F}_2$, that $V \leq \widehat{G}$ is finite and
$A \subset G$. Suppose that $\epsilon \in (0,1]$. Then either
$\|\chi_A\|_{A(G)} \geq \epsilon^{-1}$ or there is a superspace
$V'$ of $V$ with $\dim V' = \dim V+1$ and for which
\begin{equation*}
\sum_{\gamma \in V'}{|\widehat{\chi_A}(\gamma)|} \geq
\frac{\epsilon\|f_V\|_1}{2} + \sum_{ \gamma \in
V}{|\widehat{\chi_A}(\gamma)|}.
\end{equation*}
\end{lemma}
Iterating this lemma leads to the following proposition.
\begin{proposition}
Suppose that $G$ is a compact vector space over $\mathbb{F}_2$ and
$A \subset G$ is such that for all $V \leq \widehat{G}$ with $|V|
\leq M$ we have $\|f_{V}\|_1 \gg 1$. Then
\begin{equation*}
\|\chi_A\|_{A(G)} \gg \sqrt{\log M}.
\end{equation*}
\end{proposition}
We omit the proof (it is not difficult and all the ideas are
contained in the proof of Proposition \ref{finintprop}) since the
hypotheses the proposition assumes on $A$ are prohibitively
strong; nevertheless we can make use of these ideas.

\subsection{A stronger iteration lemma}

The main weakness of the above approach is that each time we apply
the weak iteration lemma to find characters supporting more
$\ell^1$-mass of $\widehat{\chi_A}$ (assuming we are not in the
case when $\|\chi_A\|_{A(G)}$ is automatically large) we do not
find very much $\ell^1$-mass, in fact we find mass in proportion
to $\|f_V\|_1$ which consequently has to be assumed large. We can
improve this by adding to $V$ not just one character at which
$\widehat{f_V}$ is large but all such characters. This idea
wouldn't work but for two essential facts.
\begin{enumerate}
\item There are a lot of characters at which $\widehat{f_V}$ is
large, in that the characters at which $\widehat{f_V}$ is large
actually support a large amount of the sum $\langle
\widehat{\chi_A},\widehat{f_V}\rangle$. \item There is a result
due to Chang which implies that the characters at which
$\widehat{f_V}$ is large are contained in a space of relatively
small dimension.
\end{enumerate}
We record Chang's result, \cite{MCC}, now. (In fact we record
something slightly different from Chang's result. For more details
see the appendix.)
\begin{theorem}\label{ChangGreen}\emph{(Chang's theorem)}
Suppose that $G$ is a compact vector space over $\mathbb{F}_2$, $f
\in L^2(G)$ and $\epsilon \in (0,1]$. Then there is a subspace $W$
of $\widehat{G}$ such that $\{\gamma:|\widehat{f}(\gamma)| \geq
\epsilon \|f\|_1 \} \subset W$ and
\begin{equation*}
\dim W \leq e\epsilon^{-2}\max\{\log
(\|f\|_2^{-2}\|f\|_1^{-2}),1\}.
\end{equation*}
\end{theorem}
\noindent We are in a position to show:
\begin{lemma}\label{finitlem}
\emph{(Iteration lemma)} Suppose that $G$ is a compact vector
space over $\mathbb{F}_2$, that $V \leq \widehat{G}$ is finite, $A
\subset G$ has $\|f_{V}\|_1>0$ and $\chi_A \in A(G)$. Then there
is a non-negative integer $s$ and a superspace $V'$ of $V$ such
that
\begin{equation*}
\sum_{\gamma \in V'}{|\widehat{\chi_A}(\gamma)|} - \sum_{ \gamma
\in V}{|\widehat{\chi_A}(\gamma)|} \gg
\left(\frac{4}{3}\right)^{s}
\end{equation*}
and
\begin{equation*}
\dim V' - \dim V \ll 4^s\log \|f_{V}\|_1^{-1}.
\end{equation*}
\end{lemma}
\begin{proof}
By Plancherel's theorem we have
\begin{equation*}
\sum_{\gamma
\in\wh{G}}{\wh{\chi_A}(\gamma)\overline{\wh{f_V}(\gamma)}} =
\langle \chi_A,f_V\rangle = \frac{1}{2}\|f_V\|_1,
\end{equation*}
where the second equality is Lemma \ref{t_Alem}. To make use of
this we apply the triangle inequality to the left hand side and
get the driving inequality of the lemma
\begin{equation}\label{leverage}
\frac{1}{2}\|f_V\|_1 \leq \sum_{\gamma \in
\wh{G}}{|\wh{\chi_A}(\gamma)||\wh{f_V}(\gamma)|}.
\end{equation}

Write $\mathcal{L}$ for the set of characters at which $\wh{f_V}$
is non-zero. Partition $\mathcal{L}$ by a dyadic decomposition of
the range of values of $|\wh{f_V}|$. Specifically, for each
non-negative integer $s$, let
\begin{equation*}
\Gamma_s:=\{\gamma \in \wh{G}:2^{-s}\|f_V\|_1 \geq
|\wh{f_V}(\gamma)| > 2^{-(s+1)}\|f_V\|_1\}.
\end{equation*}
For all characters $\gamma$ we have $|\wh{f_V}(\gamma)| \leq
\|f_V\|_1$ and if $\gamma \in \mathcal{L}$ then
$|\wh{f_V}(\gamma)| >0$ so certainly the $\Gamma_s$s cover
$\mathcal{L}$; they are clearly disjoint and hence form a
partition of $\mathcal{L}$. Write $L_s$ for the $\ell^1$-norm of
$\wh{\chi_A}$ supported on $\Gamma_s$:
\begin{equation*}
L_s:=\sum_{\gamma \in \Gamma_s}{|\wh{\chi_A}(\gamma)|}.
\end{equation*}
The right hand side of (\ref{leverage}) can now be rewritten using
these definitions:
\begin{eqnarray*}
\sum_{\gamma \in \wh{G}}{|\wh{\chi_A}(\gamma)||\wh{f_V}(\gamma)|}&
= & \sum_{\gamma \in
\mathcal{L}}{|\wh{\chi_A}(\gamma)||\wh{f_V}(\gamma)|} \textrm{ by
the definition of }\mathcal{L}\\& =
&\sum_{s=0}^\infty{\sum_{\gamma \in
\Gamma_s}{|\wh{\chi_A}(\gamma)||\wh{f_V}(\gamma)|}}\textrm{ since }\{\Gamma_s\}_{s \geq 0}\textrm{ is a partition of }\mathcal{L},\\
& \leq & \sum_{s=0}^\infty{\sum_{\gamma \in
\Gamma_s}{|\wh{\chi_A}(\gamma)|.2^{-s}\|f_V\|_1}}\textrm{ by the definition of }\Gamma_s,\\
& = & \sum_{s=0}^\infty{L_s2^{-s}\|f_V\|_1}\textrm{ by the
definition of }L_s.
\end{eqnarray*}
Combining this with (\ref{leverage}) and dividing by $\|f_V\|_1$
(which is possible since $\|f_V\|_1>0$) we get
\begin{equation}\label{averaging}
\frac{1}{2} \leq \sum_{s=0}^\infty{2^{-s}L_s}.
\end{equation}
Now, if for every non-negative integer $s$ we have
\begin{equation*}
L_s < \frac{1}{6}\left(\frac{4}{3}\right)^s,
\end{equation*}
then
\begin{equation*}
\sum_{s=0}^\infty{2^{-s}L_s} <
\sum_{s=0}^\infty{2^{-s}\frac{1}{6}\left(\frac{4}{3}\right)^s} =
\frac{1}{6}\sum_{s=0}^\infty{\left(\frac{2}{3}\right)^{s}} =
\frac{1}{2},
\end{equation*}
which contradicts (\ref{averaging}). Hence there is a non-negative
integer $s$ such that
\begin{equation*}
L_s \geq \frac{1}{6}\left(\frac{4}{3}\right)^s.
\end{equation*}

Chang's theorem gives a space $W$ for which
\begin{equation*}
\Gamma_s \subset \{\gamma \in \wh{G}:|\wh{f_V}(\gamma)| \geq
2^{-(s+1)}\|f_V\|_1\} \subset W
\end{equation*}
and
\begin{equation*}
\dim W \leq 4e2^{2s}\max\{\log
(\|f_{V}\|_2^2\|f_{V}\|_1^{-2}),1\}.
\end{equation*}
To tidy this up we note that $f_V=\chi_A - \chi_A \ast
\mu_{V^\perp}$ a.e. and $\chi_A(x),\chi_A \ast \mu_{V^\perp}(x)
\in [0,1]$, so $f_V(x) \in [-1,1]$ for a.e. $x \in G$ and hence
$\|f_V\|_2 \leq 1$, from which it follows that
\begin{equation*}
\dim W \ll 4^{s}\log \|f_V\|_1^{-1}.
\end{equation*}
Let $V'$ be the space generated by $V$ and $W$. Then
\begin{equation*}
\dim V' - \dim V \ll 4^{s} \log \|f_{V}\|_1^{-1}.
\end{equation*}

Finally we note that $\Gamma_s \cap V = \emptyset$ since
$\widehat{f_{V}}(\gamma)=0$ if $\gamma \in V$ (recall
$\widehat{f_{V}}$ from (\ref{tacomp})) and $|\wh{f_V}(\gamma)| >
2^{s+1}\|f_V\|_1\geq 0$ if $\gamma \in \Gamma_s$. Hence
\begin{equation*}
\sum_{\gamma \in V'}{|\widehat{\chi_A}(\gamma)|} \geq \sum_{\gamma
\in \Gamma_s}{|\widehat{\chi_A}(\gamma)|} + \sum_{\gamma \in
V}{|\widehat{\chi_A}(\gamma)|}  \geq
\frac{1}{6}\left(\frac{4}{3}\right)^s +\sum_{\gamma \in
V}{|\widehat{\chi_A}(\gamma)|} .
\end{equation*}
This gives the result.
\end{proof}
\noindent By iterating this lemma we prove the following result.
\begin{proposition}\label{finintprop}
Suppose that $G$ is a compact vector space over $\mathbb{F}_2$ and
$A \subset G$ is such that for all $V \leq \widehat{G}$ with $|V|
\leq M$ we have $\log \|f_{V}\|_1^{-1} \ll \log |V|$. Then
\begin{equation*}
\|\chi_A\|_{A(G)} \gg \log \log M.
\end{equation*}
\end{proposition}
\begin{proof}
Fix $\epsilon \in (0,1]$ to be optimized later. We construct a
sequence $V_0 \leq V_1 \leq ... \leq \widehat{G}$ iteratively,
writing $d_i:=\dim V_i$ and
\begin{equation*}
L_i=\sum_{\gamma \in V_i}{|\widehat{\chi_A}(\gamma)|}.
\end{equation*}
We start the construction by letting $V_0:=\{0_{\widehat{G}}\}$.
Suppose that we are given $V_k$. If $|V_k| \leq M$ then apply the
iteration lemma to $V_k$ and $A$ to get an integer $s_{k+1}$ and
vector space $V_{k+1}$ with
\begin{equation}\label{lemma conditions}
d_{k+1} - d_k \ll 4^{s_{k+1}}\log \|f_{V_k}\|_1^{-1} \textrm{ and
} L_{k+1} - L_k \gg \left(\frac{4}{3}\right)^{s_{k+1}}.
\end{equation}
First we note that the iteration terminates since certainly $L_k
\gg k$, but also $L_k \leq \|\chi_A\|_{A(G)}< \infty$.

Since $\log \|f_{V_k}\|_1^{-1} \ll \log |V_k| \ll d_k$ it follows
from (\ref{lemma conditions}) that
\begin{equation}\label{dimension bound}
d_{k+1} \ll 4^{s_{k+1}}d_k,
\end{equation}
from which, in turn, we get
\begin{equation}\label{l1 norm bound} L_k \gg
\sum_{l=0}^k{\left(\frac{4}{3}\right)^{s_l}} \gg \sum_{l
=0}^k{s_l} \gg \log d_k.
\end{equation}

Let $K$ be the stage of the iteration at which it terminates i.e.
$|V_K| > M$. We have two possibilities.
\begin{enumerate}
\item $d_{K-1}\equiv\log_2|V_{K-1}| \leq \sqrt{\log M}$: in which
case $d_{K} \geq \sqrt{\log M}.d_{K-1}$. (\ref{dimension bound})
then tells us that $4^{s_{K}}\gg \sqrt{\log M}$. However the first
inequality in (\ref{l1 norm bound}) tells us that
$\|\chi_A\|_{A(G)} \geq L_K \gg (4/3)^{S_K}$ and so certainly
$\|\chi_A\|_{A(G)} \gg \log \log M$. \item Alternatively
$d_{K-1}\equiv\log_2|V_{K-1}| \geq \sqrt{\log M}$: in which case
by (\ref{l1 norm bound}) we have $L_{K-1} \gg \log d_{K-1} \gg
\log \log M$ and so certainly $\|\chi_A\|_{A(G)} \gg \log \log M$.
\end{enumerate}
In either case the proof is complete.
\end{proof}

\section{Physical space estimates}\label{finphysical}

To realize the hypothesis of Proposition \ref{finintprop}
regarding $f_{V}$ as a density condition we have the following
lemma:
\begin{lemma}\label{lem1}
Suppose that $G$ is a compact vector space over $\mathbb{F}_2$,
that $V \leq \widehat{G}$ is finite and $A \subset G$ has density
$\alpha$. Then
\begin{equation*}
\|f_{V}\|_1 \geq 2|V|^{-1}\{\alpha |V|\}(1-\{\alpha |V|\}).
\end{equation*}
\end{lemma}
We need the following technical lemma:
\begin{lemma}\label{techlem}
Let $\delta_1,...,\delta_m \in [0,1]$ and put
$\gamma=\{\sum_{i=1}^m{\delta_i}\}$. Then
\begin{equation}\label{crit}
\sum_{i=1}^m{\delta_i-\delta_i^2} \geq \gamma(1-\gamma).
\end{equation}
\end{lemma}
\begin{proof} We may assume that $0 < \gamma < 1$. Suppose that we have $i \neq j$ such that $0<
\delta_i,\delta_j < 1$. Put $\delta=\delta_i + \delta_j \leq 2$
and we have two cases:
\begin{enumerate}
\item $ \delta \leq 1$: In this case we may replace $\delta_i$ and
$\delta_j$ by $\delta$ and 0. This preserves $\gamma$ and since
\begin{equation*}
\delta_i - \delta_i^2+ \delta_j  - \delta_j^2 \geq
(\delta_i+\delta_j) - (\delta_i+\delta_j)^2  + 0- 0^2,
\end{equation*}
it does not increase the sum in (\ref{crit}). \item $ 2\geq \delta
> 1$: In this case we may replace $\delta_i$ and $\delta_j$ by 1
and $\delta-1$. This preserves $\gamma$ and since
\begin{eqnarray*}
&(\delta_i -1)(\delta_j-1)& \geq 0 \\
\Rightarrow & 0&\geq -2\delta_i\delta_j + 2(\delta_i + \delta_j) -2\\
\Rightarrow & \delta_i - \delta_i^2+ \delta_j  - \delta_j^2 &\geq
\delta_i - \delta_i^2+ \delta_j  - \delta_j^2-2\delta_i\delta_j + 2(\delta_i + \delta_j) -2\\
\Rightarrow & \delta_i - \delta_i^2+ \delta_j  - \delta_j^2 &\geq
(\delta_i+ \delta_j - 1)  - (\delta_i+\delta_j -1)^2 + 1-1^2,
\end{eqnarray*}
it does not increase the sum in (\ref{crit}).
\end{enumerate}
In both cases we can reduce the number of $i$s for which
$0<\delta_i<1$ without increasing the sum in (\ref{crit}), so we
may assume that there is only one $j$ such that $0<\delta_j<1$.
Then
\begin{equation*}
\delta_j + \sum_{i\neq j}{\delta_i} = \gamma + \lfloor
\sum_{i=1}^m{\delta_i} \rfloor \Rightarrow \delta_j - \gamma =
\lfloor \sum_{i=1}^m{\delta_i} \rfloor -\sum_{i\neq j}{\delta_i},
\end{equation*}
but the right hand side is an integer and $-1 < \delta_j -\gamma <
1$ so $\delta_j = \gamma$ and (\ref{crit}) follows.
\end{proof}

\begin{proof}[of Lemma~{\rm\ref{lem1}}] Lemma \ref{t_Alem} states that $\|f_{V}\|_1=2\langle
\chi_A,f_{V}\rangle$ so
\begin{eqnarray*}
\|f_{V}\|_1 & = & 2\int{\chi_A(\chi_A - \chi_A \ast
\mu_{V^\perp})d\mu_{G}}\\ & = & 2\int_{x\in G}{\int{\chi_A(\chi_A
- \chi_A \ast \mu_{V^\perp})d\mu_{x+V^\perp}}d\mu_{G}(x)}.
\end{eqnarray*}
(This is just conditional expectation.) $\chi_A \ast
\mu_{V^\perp}$ is constant on cosets of $V^\perp$ and
$\chi_A^2\equiv\chi_A$ so that
\begin{eqnarray*}
\|f_{V}\|_1 & = & 2\int_{x\in G}{\int{\chi_A d\mu_{x+V^\perp}}(1 -
\chi_A \ast
\mu_{V^\perp}(x))d\mu_{G}(x)}\\
& = & 2\int_{x\in G}{\chi_A \ast \mu_{V^\perp}(x)(1 - \chi_A \ast
\mu_{V^\perp}(x))d\mu_{G}(x)}.
\end{eqnarray*}
There are $|V|$ cosets of $V^\perp$ in $G$, and $\chi_A \ast
\mu_{V^\perp}$ is constant on cosets of $V^\perp$ so this integral
is really a finite sum with $|V|$ terms in it. Let $\mathcal{C}$
be a set of coset representatives for $V^\perp$ in $G$ then
$|\mathcal{C}|=|V|$ and
\begin{equation*}
\|f_{V}\|_1 = \frac{2}{|\mathcal{C}|} \sum_{x' \in
\mathcal{C}}{\chi_A \ast \mu_{V^\perp}(x')(1 - \chi_A \ast
\mu_{V^\perp}(x'))}.
\end{equation*}
We can now apply Lemma \ref{techlem} to the $\chi_A \ast
\mu_{V^\perp}(x')$s with $m=|\mathcal{C}|$. This gives
\begin{equation*}
\|f_{V}\|_1 \geq \frac{2}{|\mathcal{C}|}\beta(1-\beta)=
\frac{2}{|V|}\beta(1-\beta)
\end{equation*}
where
\begin{equation*}
\beta=\left\{\sum_{x' \in \mathcal{C}}{\chi_A \ast
\mu_{V^\perp}(x')}\right\} =\left\{|\mathcal{C}|\int_{x \in
G}{\chi_A \ast
\mu_{V^\perp}(x)d\mu_G(x)}\right\}=\left\{|V|\alpha\right\}.
\end{equation*}
\end{proof}
Nothing better than Lemma \ref{lem1} can be true: Let $A$ be the
union of $\lfloor \alpha|V| \rfloor $ cosets of $V^\perp$ and a
subset of a coset of $V^\perp$ of relative density
$\{\alpha|V|\}$. Equality is attained in Lemma \ref{lem1} for this
set.

\section{The result, remarks and examples}\label{finrem}

As an easy corollary of Proposition \ref{finintprop} and Lemma
\ref{lem1} we have:
\begin{theorem}\label{f2ndodge}
Suppose that $G$ is a compact vector space over $\mathbb{F}_2$.
Suppose that $A \subset G$ has density $\alpha$ and for all $V
\leq \widehat{G}$ with $|V|\leq M$ we have
$\{\alpha|V|\}(1-\{\alpha|V|\}) \gg |V|^{-1}$. Then
\begin{equation*}
\|\chi_A\|_{A(G)} \gg \log \log M.
\end{equation*}
\end{theorem}
Theorem \ref{maintheorem} is simply a weaker version of this
result.

If we were interested we could read the explicit dependencies
between the implied constants in Theorem \ref{f2ndodge} out of the
proof. Since the only real importance of this result is in its
corollary, Theorem \ref{maintheorem}, where even the $M$
dependence is almost certainly not best possible, this is probably
of little interest.

There are clear similarities between the methods of this paper and
those employed by Green and Konyagin in \cite{BJGSVK}, however the
most striking ones seem to be between this work and the work of
Bourgain in \cite{JBIS}. In particular a slight variation on the
calculation in Lemma \ref{t_Alem} is in his work and he proves a
result using Beckner's inequality (which is essentially equivalent
to Chang's theorem) which shows that if $A \subset \mathbb{F}_2^n$
has density $\alpha$ with $\alpha(1-\alpha) \gg 1$ then either
$\widehat{\chi_A}$ is large at a non-trivial character or there is
significant $\ell^2$-mass in the tail of the Fourier transform.

Our reason for stating Theorem \ref{f2ndodge} at all is that it is
sharp up to the constant and hence demonstrates a limitation of
our method. If we let
\begin{equation*}
\alpha=\frac{1}{2^{2^0}}+\frac{1}{2^{2^1}}+...+\frac{1}{2^{2^{k-1}}},
\end{equation*}
then we showed in \S\ref{egs} that there is a set $A$ of density
$\alpha$ with $\|\chi_A\|_{A(G)} \leq k$. However $A$ also
satisfies the hypotheses of Theorem \ref{f2ndodge} with
$M=2^{2^{k-1}}-1$: If $V \leq \wh{G}$ has $|V| \leq M$ then
$|V|=2^{d}$ for some $d <2^{k-1}$,
\begin{eqnarray*}
\{\alpha |V|\} = \sum_{\min\{0,\log_2{d}\} <m \leq
k-1}{2^d.2^{-2^m}}& \leq & \sum_{\min\{0,\log_2{d}\} <m \leq
k-1}{2^{-2^{m-1}}}\\ & \leq & \sum_{m=0}^\infty{2^{-2^{m}}} \leq
\frac{7}{8},
\end{eqnarray*}
and
\begin{equation*}
\{\alpha |V|\} = \sum_{\min\{0,\log_2{d}\} <m \leq
k-1}{2^d.2^{-2^m}} \geq 2^d.2^{-2^{\lfloor \log_2{d} \rfloor +1}}
\geq 2^{-d} = |V|^{-1}.
\end{equation*}
Hence
\begin{equation*}
\{\alpha|V|\}(1-\{\alpha|V|\}) \gg |V|^{-1}.
\end{equation*}
Theorem \ref{f2ndodge} then tells us that $\|\chi_A\|_{A(G)} \gg k
$.

\section*{Acknowledgements}
I should like to thank Tim Gowers for supervision and for
reviewing a number of drafts of this paper, Ben Green for many
valuable conversations and comments, and for pointing out
\cite{JBIS}, and Ben Green and Sergei Konyagin for making the
preprint \cite{BJGSVK} available.

\appendix

\section{Chang's theorem in compact vector spaces over
$\mathbb{F}_2$}

Green was the first to observe that the original proof of Chang's
theorem in \cite{MCC} can easily be adapted to prove a result not
just about characteristic functions of sets, but about arbitrary
elements of $L^2(G)$. It is this result which we recorded as
Chang's theorem (Theorem \ref{ChangGreen}) above. We take the
opportunity to provide a proof of this using Beckner's inequality
which as Green has noted in \cite{BJGRKP}, can be used in place of
Rudin's inequality when we are in the setting of compact vector
spaces over $\mathbb{F}_2$.

Suppose that $G$ is a compact vector space over $\mathbb{F}_2$. If
$\Lambda$ is a finite set of characters on $G$ and $\eta \in
[-1,1]$ then we can define the Riesz product
\begin{equation*}
p_\eta:=\prod_{\lambda \in \Lambda}{(1+\eta \lambda)}.
\end{equation*}
Every term in this product is real and non-negative since $\eta
\in [-1,1]$, so $p_\eta$ is real and non-negative.

Now suppose that $\Lambda$ is linearly independent. Then
$\widehat{p_\eta}(0_{\widehat{G}})=1$ and hence $\|p_\eta\|_1 =
\widehat{p_\eta}(0_{\widehat{G}})=1$ by non-negativity of
$p_\eta$. It follows by Young's inequality that we can define the
operator
\begin{eqnarray*}
T_\eta:L^2(G) & \rightarrow & L^2(G)\\ f & \mapsto & f \ast
p_\eta,
\end{eqnarray*}
and moreover $\|T_\eta f\|_2 \leq \|f\|_2\|p_\eta\|_1 = \|f\|_2$.
In fact there is a stronger inequality.
\begin{theorem}\emph{(Beckner's inequality)}
Suppose that $G$ is a compact vector space over $\mathbb{F}_2$ and
$\Lambda$ is a finite linearly independent set of characters on
$G$. Suppose that $\eta \in [-1,1]$. Then the operator $T_\eta$
defined above has $\|T_\eta f\|_2 \leq \|f\|_{1+\eta^2}$.
\end{theorem}
Having stated Beckner's inequality we are in a position to proceed
with the proof of Chang's theorem. Chang noted the following
simple fact regarding linearly independent sets.
\begin{lemma}
Suppose that $V$ is a vector space, $\Gamma$ is a subset of $V$
and $\Lambda$ is a maximal linearly independent subset of
$\Gamma$. Then $\Gamma$ is contained in the subspace generated by
$\Lambda$.
\end{lemma}
Theorem \ref{ChangGreen} now follows from this and the following
result.
\begin{proposition}
Suppose that $G$ is a compact vector space over $\mathbb{F}_2$, $f
\in L^2(G)$ and $\epsilon \in (0,1]$. Suppose that $\Lambda$ is a
linearly independent subset of
$\Gamma:=\{\gamma:|\widehat{f}(\gamma)| \geq \epsilon \|f\|_1 \}$.
Then
\begin{equation*}
|\Lambda| \leq e\epsilon^{-2}\max\{\log
\left(\|f\|_2^{-2}\|f\|_1^{-2}\right),1\}
\end{equation*}
\end{proposition}
\begin{proof}
It certainly suffices to prove the result for $\Lambda$ finite. In
this case the operator $T_\eta$ is defined and we can apply
Beckner's inequality.
\begin{equation*}
\sum_{\lambda \in \Lambda}{|\eta\widehat{f}(\lambda)|^2} \leq
\sum_{\gamma \in \wh{G}}{|\widehat{T_\eta f}(\gamma)|^2} =
\|T_\eta f\|_2 \leq \|f\|_{1+\eta^2}.
\end{equation*}
The equality in the middle is by Parseval's theorem. Since
$\Lambda \subset \Gamma$ it follows that
\begin{equation*}
(\eta\epsilon\|f\|_1)^2|\Lambda| \leq \sum_{\lambda \in
\Lambda}{|\eta\widehat{f}(\lambda)|^2} \leq \|f\|_{1+\eta^2}.
\end{equation*}
After some manipulation and using the log-convexity of $\|.\|_p$
we get
\begin{equation*}
|\Lambda| \leq \eta^{-2}\epsilon^{-2}\|f\|_1^{-2}
\|f\|_{1+\eta^2}^2 \leq
\epsilon^{-2}\eta^{-2}\left(\frac{\|f\|_2^2}{\|f\|_1^2}\right)^{\eta^2}.
\end{equation*}
Finally we optimize by putting $\eta^{-2}=\max\{\log
(\|f\|_2^2\|f\|_1^{-2}),1\}$ to get the result.
\end{proof}

\bibliographystyle{alpha}

\bibliography{master}

\end{document}